\documentclass[final,onefignum,onetabnum]{siamart251104}
\usepackage{xcolor}
\usepackage{lipsum}
\usepackage{amsfonts}
\usepackage{graphicx}
\usepackage{epstopdf}
\usepackage{algorithmic}
\usepackage{float, soul}
\usepackage{arydshln}
\usepackage{amssymb}
\ifpdf
  \DeclareGraphicsExtensions{.eps,.pdf,.png,.jpg}
\else
  \DeclareGraphicsExtensions{.eps}
\fi

\newsiamremark{remark}{Remark}
\newsiamremark{assumption}{Assumption}
\crefname{assumption}{Assumption}{Assumptions}
\newsiamthm{claim}{Claim}
\newsiamremark{fact}{Fact}
\crefname{fact}{Fact}{Facts}

\headers{Symmetry-reduced HiSD}{Y. Lin, J. Yin, and L. Zhang}

\title{Symmetry-reduced High-index Saddle Dynamics for Saddle Point Search under Continuous Symmetries\thanks{Submitted to the editors DATE.
\funding{This work was supported by the National Natural Science Foundation of China (No. 12225102, T2321001, 12288101, 12601819).}}
}

\author{Yuanshen Lin\thanks{School of Mathematical Sciences, 
Peking University, Beijing 100871, China. 
(\email{yslin25@stu.pku.edu.cn}).}
\and 
Jianyuan Yin\thanks{School of Mathematical Sciences, 
Laboratory of Mathematics and Complex Systems, Ministry of Education, 
Beijing Normal University, Beijing 100875, China.
(\email{jyyin@bnu.edu.cn}).}
\and 
Lei Zhang\thanks{\raggedright Corresponding author. School of Mathematical Sciences, Beijing International Center for Mathematical Research, Center for Quantitative Biology, Center for Machine Learning Research, Institute for Artificial Intelligence, Peking University, Beijing, China. (\email{zhangl@math.pku.edu.cn}).}
}
\usepackage{amsopn}

\ifpdf
\hypersetup{
  pdftitle={Symmetry-reduced HiSD},
  pdfauthor={Lin et al.}
}
\fi

\begin{document}

\maketitle

\begin{abstract}
For energy invariant under continuous symmetries, critical points may form critical manifolds, leading to degeneracy of the Hessian. 
In such cases, the ascent directions in high-index saddle dynamics (HiSD) may include symmetry directions, making it difficult to locate saddle points of a prescribed Morse index. 
To address this issue, we propose symmetry-reduced HiSD (SR-HiSD) by formulating the quotient HiSD on the quotient manifold and horizontally lifting the dynamics to the ambient Hilbert space. 
Under suitable regularity and Morse--Bott assumptions, we establish local exponential stability of the quotient HiSD equilibria associated with the prescribed index and instability of equilibria associated with critical points of different indices. We further prove local convergence of the SR-HiSD position trajectory to a single point on the target critical orbit. Numerical algorithms with direction updates restricted to the horizontal space are developed. Numerical experiments on Morse clusters, the Onsager model, and the Landau--Brazovskii model demonstrate reliable searches for saddle points of the prescribed index and the construction of solution landscapes, together with improved computational efficiency in large-scale problems.
\end{abstract}



\section{Introduction}
\label{sec:1}
Locating saddle points is important for understanding complex energy landscapes. 
Local minima represent stable or metastable states, whereas index-1 saddle points characterize transition states between them.
Saddle-point computations play an important role in the study of phase transitions \cite{cheng2010nucleation,samanta2014microscopic}, chemical reactions \cite{heyden2005efficient, jacobson2017automated} and neural-network training \cite{dauphin2014identifying,zhang2021embedding}. 
However, saddle points are unstable and can be difficult to compute, especially in high-dimensional systems where transition mechanisms are unclear. 
These challenges motivate the development of efficient and reliable numerical methods for saddle-point computations.

Path-finding and surface-walking methods are two widely used classes of saddle-search methods. 
Path-finding methods, such as the nudged elastic band method \cite{henkelman2000improved} and the string method \cite{weinan2002string}, compute minimum energy paths connecting given initial and final states. 
Surface-walking methods use local derivative information to search for saddle points from an initial guess, while they do not require prior knowledge of the final state of a transition. 
Representative methods include gentlest ascent dynamics \cite{weinan2011gentlest}, dimer-type methods \cite{henkelman1999dimer}, the iterative minimization formulation \cite{gao2015iterative}, eigenvector-following methods \cite{baker1986algorithm}, gradient extremals \cite{jorgensen1988gradient} and Newton trajectories \cite{quapp1998searching}. 
Within the surface-walking framework, high-index saddle dynamics (HiSD) was developed to locate saddle points of a prescribed Morse index $k$ \cite{yin2019high}, by combining gradient ascent in an evolving $k$-dimensional subspace with gradient descent in its orthogonal complement. 
Starting from known critical points, upward and downward searches based on HiSD can locate critical points with different indices and construct solution landscapes that describe their connections \cite{yin2020construction}. 
HiSD and its extensions have been applied to liquid crystals \cite{wang2021modelling}, quasicrystals \cite{yin2021transition}, and Bose--Einstein condensates \cite{yin2024revealing}.

HiSD was originally formulated for nondegenerate saddle points, for which local stability and convergence results have been established \cite{luo2022convergence,zhang2022sinum}.
A critical point is nondegenerate if its Hessian has a bounded inverse \cite{milnor1963morse}, and such a nondegenerate critical point is isolated. 
This assumption naturally excludes critical points on positive-dimensional critical manifolds, which arise in many applications \cite{yin2022solution,yin2024revealing}.
In many variational models, the energy is invariant under continuous symmetries described by Lie group actions, such as translations, rotations, or phase shifts.  
Consequently, each physical state is represented by an entire orbit of symmetry-related points in the ambient space. 
At a critical point, the tangent space to the symmetry orbit lies in the kernel of the Hessian, so the point is degenerate whenever the orbit has positive dimension. 
These symmetry directions in the orbit tangent space describe motion along the orbit between symmetry-equivalent representatives of the same physical state \cite{golubitsky2012singularities,field2007dynamics,marsden1999introduction}. 
When the Hessian is nondegenerate on the normal space to a critical orbit, this orbit forms a Morse--Bott critical manifold \cite{bott1954nondegenerate,banyaga2010morse}. 
Despite this degeneracy, HiSD has been applied successfully to compute degenerate saddle points in many physical systems \cite{yin2022solution,yin2024revealing}.
Under suitable Morse--Bott and regularity assumptions, Luo et al. established local convergence of HiSD to degenerate saddle points on critical manifolds, while the prescribed HiSD index may not be the Morse index of the converged saddle point because the ascent subspace may include symmetry directions in addition to unstable directions \cite{luo2026convergence}. 
Recently, by explicitly computing the Hessian nullspace, methods that constrain the ascent directions to be orthogonal to this nullspace have also been developed, including the nullspace-preserving saddle search (NPSS) method for ordered phase transitions \cite{cui2025efficient} and its high-index extensions \cite{jiang2025nullspace}.
These methods construct ascent directions orthogonal to a reference nullspace which is updated between search segments.

Away from critical points, however, symmetry directions need not lie in the Hessian kernel, even though the energy remains constant along each orbit. 
This distinction motivates us to identify the symmetry directions directly from the group action and treat them separately throughout the saddle-point search.
For energy landscapes with continuous symmetries, a natural approach is to formulate the saddle-point problem on the quotient space by treating each symmetry orbit as a single object.
Under suitable assumptions on the smooth isometric group action, the quotient space is a smooth Riemannian manifold, and the quotient map is a Riemannian submersion \cite{palais1961existence,lee2013smooth,kobayashi1996foundations,o1966fundamental}. 
Tangent vectors on the quotient are represented by directions orthogonal to the symmetry orbits in the ambient space, and saddle-point searches can therefore be developed using Hessian on the quotient. 
The Morse index on the quotient is determined by the number of negative eigenvalues of this reduced Hessian.

In this paper, we reduce the continuous symmetries and formulate HiSD directly on the quotient manifold. 
The quotient dynamics are defined using the Riemannian gradient and Hessian together with an evolving orthonormal $k$-frame. 
We then derive symmetry-reduced HiSD (SR-HiSD) in the ambient space as a horizontal lift of the quotient dynamics.
Under suitable Morse--Bott and regularity assumptions, we establish local exponential stability of the quotient dynamics at an index-$k$ saddle point.
We further show that nearby ambient trajectories converge to individual index-$k$ saddle points on the corresponding critical orbit.
We also construct discrete SR-HiSD algorithms that update the state using a projected Euler step and approximate the unstable frame using a projected eigensolver.
Numerical experiments involving rigid-motion, rotational, and translational symmetries demonstrate the effectiveness of SR-HiSD in locating saddle points of the prescribed Morse index.

The remainder of the paper is organized as follows. 
\Cref{sec:main_results_q} presents the quotient formulation, derives the quotient HiSD and SR-HiSD method, and establishes their equivalence and local convergence properties.
\Cref{sec:algorithm_q} describes the discrete SR-HiSD algorithms and step-size selection.
\Cref{sec:experiments_q} presents numerical experiments. 
Conclusions and remarks are presented in \Cref{sec:5}.

\section{SR-HiSD}\label{sec:main_results_q}
In this section, we develop the HiSD method under the quotient-manifold framework. 
The main idea is to remove symmetry-induced degeneracy at the level of the state space itself. 
Instead of dealing with a degenerate critical point in the ambient Hilbert space, we regard it as a nondegenerate critical point on the quotient manifold, together with its horizontal lift back to the ambient space. 
We begin with the quotient-manifold setting, then formulate the intrinsic HiSD in the quotient space, and finally return to the ambient representation of HiSD and its convergence properties.

\subsection{Quotient Manifold}
Let $H$ be a real Hilbert space with inner product $\langle \cdot,\cdot \rangle$, and $G$ be a finite-dimensional Lie group acting smoothly and isometrically on $H$. 
Let $\Omega\subset H$ be a $G$-invariant open set and assume the group action on $\Omega$ is free and proper.
For $x\in\Omega$, we denote its orbit by $G\cdot x = \{g\cdot x: g\in G\} \subset \Omega$, which is a smooth manifold of dimension $\dim G$. 
Then we denote the quotient space by $\mathcal M := \Omega/G$ and the quotient map by $\pi :\Omega \to \mathcal M$. 
For $x\in\Omega$, the vertical space $\mathcal V_x$ is defined as $\mathcal V_x := \ker D\pi_x = T_x(G\cdot x)$.
Then the tangent space $T_x H$ admits an orthogonal decomposition $T_x H = \mathcal V_x \oplus \mathcal H_x$, where $\mathcal H_x:=\mathcal V_x^\perp$ is the horizontal space. 
For proper actions in Hilbert settings, the closedness of orbit tangent spaces is a useful technical point in making this orthogonal decomposition well behaved \cite{jotz2008closedness}. 
We denote the orthogonal projection onto $\mathcal H_x$ by $P_x: T_xH\to\mathcal H_x$. 
Since the group action is isometric, the restriction of $D\pi_x$ to $\mathcal H_x$ is an isometry onto $T_{[x]}\mathcal M$, so we define the horizontal lift by $(\cdot)^\sharp: T_{[x]}\mathcal M \to \mathcal H_x$, i.e., for $\eta\in T_{[x]}\mathcal M$, $\eta^\sharp\in\mathcal H_x$ is the unique horizontal vector such that $D\pi_x \eta^\sharp = \eta$ \cite{o1966fundamental}.

Throughout this paper, we work in a regular neighborhood of a principal orbit under the following assumption. 
\begin{assumption}[Regular quotient setting]
\label{hyp:regular_quotient}
The action of $G$ on $\Omega$ is smooth, free, proper, and isometric.
Consequently, the quotient space $\mathcal M=\Omega/G$ is a
smooth Hilbert manifold, and the quotient map $\pi$ is a
smooth submersion
\cite{palais1961existence,jotz2008closedness}. Moreover, the metric on $H$ induces a Riemannian metric
$\langle\cdot,\cdot\rangle_{\mathcal M}$ on $\mathcal M$
for which $\pi$ is a Riemannian submersion
\cite{kobayashi1996foundations,o1966fundamental}.
Therefore, $D\pi_{x}|_{\mathcal H_{x}}$ is an isometric isomorphism onto
$T_{[x]}\mathcal M$ and the horizontal projector $x\mapsto P_x$ is of class $C^1$ as an operator-valued map \cite{boumal2020introduction}.
\end{assumption}

Assumption \ref{hyp:regular_quotient} ensures the smooth quotient structure, the regularity of the horizontal projector, and the Riemannian submersion used below. 
Near singular orbits, these objects generally require a stratified-space formulation.

We consider an energy functional $E:\Omega\to\mathbb R$ which is of class $C^2$. We first recall some standard notions from infinite-dimensional Morse theory \cite{chang1993infinite}. A point $x^*\in\Omega$ is called a critical point of $E$ if its gradient vanishes, i.e., $\nabla E(x^*)=0$. 
The Morse index of a critical point $x^*$, denoted by $\operatorname{ind}_E(x^*)$, is the maximal dimension of a subspace of $H$ on which the quadratic form associated with the Hessian $\nabla^2E(x^*)$ is negative definite. Thus, if $\operatorname{ind}_E(x^*)=k$, then $x^*$ is called an index-$k$ critical point of $E$.
The nullity of $x^*$ is defined as $\dim\ker\nabla^2E(x^*)$.  
A critical point $x^*$ is termed nondegenerate if the Hessian $\nabla^2E(x^*)$ has a bounded inverse.

We now assume that the energy functional $E$ is $G$-invariant, namely,
\begin{equation}
E(g\cdot x)=E(x), \qquad \forall g\in G,\; x\in \Omega.
\end{equation}
Because $E$ is constant on each orbit, we can define a reduced energy
\begin{equation}
  \mathcal E:\mathcal M\to\mathbb R,
  \qquad
  \mathcal E([x]) := E(x).
\end{equation}
Moreover, the gradient of $E$ is horizontal (i.e., $\nabla E(x) \in \mathcal{H}_x$), and the gradient of $\mathcal E$ is characterized by
\begin{equation}\label{eq:reduced_gradient}
  \bigl(\operatorname{grad}\mathcal E([x])\bigr)^\sharp = \nabla E(x).
\end{equation}
Likewise, for any $\eta \in T_{[x]}\mathcal M$ with horizontal lift $\eta^\sharp\in\mathcal H_x$, the Hessian satisfies
\begin{equation}\label{eq:reduced_hessian}
\bigl(\operatorname{Hess}\mathcal E([x])\eta\bigr)^\sharp
= P_x\nabla^2E(x)\eta^\sharp.
\end{equation}
Equation \eqref{eq:reduced_hessian} is the standard Hessian identity for an invariant function under a Riemannian submersion; the curvature correction vanishes because $E$ is constant along the fibers \cite{o1966fundamental,abraham2012manifolds}.

We point out that a quotient viewpoint is natural for symmetry-reduced saddle search on the $G$-invariant energy $E$. 
Since $G$ is free and isometric on $\Omega$, for an index-$k$ saddle point $x^*\in\Omega$, the orbit $G\cdot x^*$ is a manifold of dimension $\dim G$, and each point on this orbit is equivalently an index-$k$ saddle point of $E$ with the same energy and local geometry. 
Because $\mathcal V_{x^*}\subset \ker \nabla^2 E(x^*)$, the Hessian $\nabla^2 E(x^*)$ has nullity at least $\dim G$. 
In other words, the critical point $x^*\in\Omega$ is typically degenerate and not isolated. 
Therefore, to formulate the nondegeneracy on the quotient, we introduce the Morse--Bott theory \cite{chang1993infinite,bott1954nondegenerate}. 
Let $\mathcal C\subset\Omega$ be a smooth critical manifold of $E$, i.e., each point on $\mathcal C$ is a critical point of $E$. 
Since $\nabla E$ vanishes on $\mathcal C$, one has 
\begin{equation}
T_x\mathcal C \subset \ker\nabla^2E(x), \qquad \forall x\in\mathcal C. 
\end{equation}
In other words, critical points on the critical manifold $\mathcal C$ are degenerate. 
The critical manifold $\mathcal C$ is termed nondegenerate in the sense of Morse--Bott if 
\begin{equation}\label{eq:morse_bott_kernel} 
\ker\nabla^2 E(x) = T_x\mathcal C, \qquad \forall x\in\mathcal C, 
\end{equation}
namely, the Hessian is nondegenerate along directions in the normal space of $\mathcal C$ \cite{bott1954nondegenerate}.

Since the action is free, the quotient removes precisely the symmetry directions
$\mathcal V_{x^*}$, with $\dim\mathcal V_{x^*}=\dim G$.
At a critical point, $\mathcal V_{x^*}\subset
\ker\nabla^2E(x^*)$, so these directions contribute
to the degeneracy of the Hessian.
By Assumption \ref{hyp:regular_quotient} and Equation \eqref{eq:reduced_hessian}, the negative
directions of $\nabla^2E(x^*)$ correspond exactly to those of
$\operatorname{Hess}\mathcal E([x^*])$.
Hence, we have
\begin{equation}\label{eq:index_correspondence}
\operatorname{ind}_E(x^*)
=
\operatorname{ind}_{\mathcal E}([x^*]).
\end{equation}

For a critical point $x^*$ of the $G$-invariant energy $E$, the manifold $\mathcal C=G\cdot x^*$ is called a critical orbit, which is represented on the quotient manifold by the critical point $[x^*]$ of $\mathcal E$. We call $G\cdot x^*$ a nondegenerate critical orbit if it is a nondegenerate critical manifold in the sense of Morse-Bott defined by Equation \eqref{eq:morse_bott_kernel}.

The Hessian restricted to the horizontal space is
\begin{equation}\label{eq:transverse_hessian}
    L_{x^*}
    :=
    P_{x^*}\nabla^2E(x^*)P_{x^*}
    \big|_{\mathcal H_{x^*}}
    :
    \mathcal H_{x^*}\to\mathcal H_{x^*}.
\end{equation}
By \eqref{eq:reduced_hessian}, $L_{x^*}$ corresponds to
$\operatorname{Hess}\mathcal E([x^*])$ through the horizontal lift.

To ensure the bounded invertibility of the Hessian required for nondegeneracy,  
we therefore impose a spectral gap of the transverse Hessian around zero.

\begin{assumption}
\label{hyp:transversal_morse}
For the nondegenerate critical orbit $G\cdot x^*$, there exist
$\gamma_-,\gamma_+>0$ such that the spectrum of $L_{x^*}$  satisfies
\begin{equation}
    \sigma(L_{x^*})
    \subset
    (-\infty,-\gamma_-]\cup[\gamma_+,\infty).
\end{equation}
\end{assumption}

Under Assumption \ref{hyp:transversal_morse}, if
$\operatorname{ind}_{\mathcal E}([x^*])=k$, we call $G\cdot x^*$ a nondegenerate index-$k$ critical orbit of $E$. 
Each point on the orbit is a degenerate index-$k$ critical point of $E$ in the ambient space $H$, where the degeneracy is generated precisely by the symmetry directions $\mathcal V_{x^*}$.

\subsection{Review of HiSD}
We first recall high-index saddle dynamics (HiSD), which searches
for an index-$k$ saddle point through a local min--max formulation
\cite{yin2019high,luo2026convergence}.
For a smooth energy $E:\mathbb R^n\to\mathbb R$, 
HiSD for an index-$k$ saddle point ($k$-HiSD) is
\begin{equation}\label{eq:classical_hisd}
\left\{
\begin{aligned}
    \frac{\mathrm dx}{\mathrm dt}
    &=
    -\left(I-2\sum_{j=1}^k v_jv_j^\top\right)\nabla E(x),
    \\
    \frac{\mathrm dv_i}{\mathrm dt}
    &=
    -\left(I-v_iv_i^\top
    -2\sum_{j=1}^{i-1}v_jv_j^\top\right)\nabla^2E(x)v_i,
    \qquad i=1,\ldots,k,
\end{aligned}
\right.
\end{equation}
with $\langle v_i(0),v_j(0)\rangle=\delta_{ij}$.

The $x$-dynamics in Equation \eqref{eq:classical_hisd} performs gradient ascent on
$\operatorname{span}\{v_1,\ldots,v_k\}$ and gradient descent along
its orthogonal complement. 
The $v_i$-dynamics minimizes the Rayleigh quotient 
$\langle\nabla^2E(x)v_i,v_i\rangle$ under the orthonormal constraint $\langle v_i,v_j\rangle=\delta_{ij}$ for $j\leq i$
to calculate the eigenvector corresponding to the $i$th smallest eigenvalue of the Hessian $\nabla^2E(x)$.
The numerical analysis and convergence results of continuous dynamics and discrete algorithms can be found in \cite{yin2019high,luo2022convergence}. 

In this $G$-invariant setting, the quotient map removes the symmetry directions while preserving the saddle index by \eqref{eq:index_correspondence}.
Under Assumption \ref{hyp:transversal_morse}, we can apply similar ideas to search for saddle points of the reduced energy $\mathcal E$ on the quotient manifold $\mathcal M$.

\subsection{HiSD on the Quotient Manifold} 
We now formulate HiSD intrinsically under the quotient setting to search for an index-$k$ saddle point on the quotient manifold $\mathcal{M}$. The position variable is represented by a point $[x]\in\mathcal M$, while the directional variables are given by an ordered orthonormal $k$-frame
$\mathbf V=(V_1,\dots,V_k)$ in $T_{[x]}\mathcal M$. We define the associated Stiefel bundle by
\begin{equation}
    \operatorname{St}_k(T\mathcal M)
    :=
    \left\{
        ([x],\mathbf V):
        \begin{array}{l}
            [x]\in\mathcal M,\quad
            \mathbf V=(V_1,\ldots,V_k)
            \in(T_{[x]}\mathcal M)^k,\\
            \langle V_i,V_j\rangle_{\mathcal M}
            =\delta_{ij},
            \quad 1\le i,j\le k
        \end{array}
    \right\}.
\end{equation}
The fiber of $\operatorname{St}_k(T\mathcal M)$ over $[x]$ is
$\operatorname{St}_k(T_{[x]}\mathcal M)$, equipped with its
standard Hilbert-manifold structure \cite{lee2013smooth,kobayashi1996foundations,harms2012geodesics}.

Following the constructions of Equation \eqref{eq:classical_hisd}, for every $([x],\mathbf V)\in \operatorname{St}_k(T\mathcal M)$, define the reflection operator on $T_{[x]}\mathcal M$ by $\mathcal A(\mathbf V): T_{[x]}\mathcal M \to T_{[x]}\mathcal M$,
\begin{equation}
  \mathcal A(\mathbf V)\eta
  := \eta - 2\sum_{j=1}^k \langle \eta, V_j\rangle_{\mathcal M} V_j.
\end{equation}
For $i=1,\dots,k$, define the deflated projection operator as $\mathcal P_i(\mathbf V): T_{[x]}\mathcal M \to T_{[x]}\mathcal M$,
\begin{equation}
  \mathcal P_i(\mathbf V)\eta
  := \eta - \langle \eta, V_i\rangle_{\mathcal M}V_i
     - 2\sum_{j=1}^{i-1}\langle \eta, V_j\rangle_{\mathcal M}V_j.
\end{equation}
Then the quotient HiSD for an index-$k$ saddle point of $\mathcal E$ is
\begin{equation}\label{eq:new_intrinsic_hisd}
\left\{\begin{aligned}
  \frac{\mathrm d[x]}{\mathrm dt}
  &= -\mathcal A(\mathbf V)\operatorname{grad}\mathcal E([x]), \\
  \frac{\mathrm DV_i}{\mathrm dt}
  &= -\mathcal P_i(\mathbf V)\bigl(\operatorname{Hess}\mathcal E([x])V_i\bigr),
  \qquad i=1,\dots,k,
\end{aligned}\right.
\end{equation}
where $\frac{\mathrm D}{\mathrm dt}$ denotes the Levi--Civita covariant derivative along the trajectory $[x](t)$.

It is straightforward to find that a stationary state of Equation \eqref{eq:new_intrinsic_hisd} is a pair $([x]^*,\mathbf V^*)\in \operatorname{St}_k(T\mathcal M)$ such that $[x]^*$ is a critical point of $\mathcal E$ and $\mathbf V^* = (\phi_1,\dots, \phi_k)$ is an ordered orthonormal eigenframe of $\operatorname{Hess}\mathcal E([x]^*)$, i.e., 
\begin{equation}
\operatorname{Hess}\mathcal E([x]^*) \phi_j =\langle \operatorname{Hess}\mathcal E([x]^*) \phi_j, \phi_j\rangle_{\mathcal M}\phi_j
\qquad
\langle \phi_i, \phi_j \rangle_{\mathcal M} = \delta_{ij}
. 
\end{equation}
In particular, because the operators $\mathcal P_i$ use sequential deflation, an arbitrary orthogonal transformation  of basis in an invariant subspace need not be a stationary frame.
Since the dynamics \eqref{eq:new_intrinsic_hisd} evolves on $\operatorname{St}_k(T\mathcal M)$, stability is understood with respect to the local bundle structure. To pass from the linear exponential stability established below to nonlinear local exponential stability, we impose the following standard regularity assumption. It is automatic in finite dimensions and standard for semilinear evolution equations
\cite{henry2006geometric,pazy2012semigroups}.

\begin{assumption}[Local semiflow regularity]
\label{hyp:local_semiflow}
In a local chart of $\operatorname{St}_k(T\mathcal M)$ near $([x]^*,\mathbf V^*)$, the quotient HiSD vector field in Equation \eqref{eq:new_intrinsic_hisd} is $C^1$ and generates a local $C^1$ semiflow. This regularity holds, for example, when the energy
is $C^3$ and the underlying geometric structures
are sufficiently smooth.

\end{assumption}

For notational convenience, we write the vector field associated
with \eqref{eq:new_intrinsic_hisd} by
\begin{equation}
\label{eq:quotient_hisd_vector_field}
    F([x],\mathbf V)
    =
    \begin{pmatrix}
        F_x\\
        F_1 \\
        \vdots \\
        F_k
    \end{pmatrix}
    :=
    \begin{pmatrix}
        -\mathcal A(\mathbf V)\operatorname{grad}\mathcal E([x])\\
        -\mathcal P_1(\mathbf V)\operatorname{Hess}\mathcal E([x])V_1\\
        \vdots\\
        -\mathcal P_k(\mathbf V)\operatorname{Hess}\mathcal E([x])V_k
    \end{pmatrix}. 
\end{equation}

We first consider the frame dynamics with the position fixed at the critical point.
\begin{lemma} 
\label{lem:frame_linearization} 
Under Assumptions \ref{hyp:regular_quotient} and \ref{hyp:transversal_morse}, let $[x]^*\in\mathcal M$ be a nondegenerate index-$k$ critical point. Assume that the $k$ negative eigenvalues of $L:=\operatorname{Hess}\mathcal E([x]^*)$ are ordered as $\lambda_1<\cdots<\lambda_k<0,$ with the corresponding orthonormal eigenframe $\mathbf V^*=(\phi_1,\dots,\phi_k)$. Then, with $[x]=[x]^*$ fixed, $\mathbf V^*$ is a linearly exponentially stable equilibrium of the frame dynamics in Equation \eqref{eq:new_intrinsic_hisd} in $\operatorname{St}_k(T_{[x]^*}\mathcal M)$. \end{lemma}

\begin{proof}
Let $W^-=\operatorname{span}\{\phi_1,\dots,\phi_k\}$,
$W^+=(W^-)^\perp$, and $T_{[x]^*}\mathcal M=W^-\oplus W^+$.
We study the linearized frame dynamics in the tangent space $T_{\mathbf V^*}\operatorname{St}_k(T_{[x]^*}\mathcal M)$, whose norm is defined as $\left\| \cdot \right\|_{V^*}$. Let $\delta\mathbf V=(\delta V_1,\dots,\delta V_k)$ be a tangent variation at $\mathbf V^*$. To obtain $\langle\phi_i,\phi_j\rangle_{\mathcal M} = \delta_{ij}$, the tangent variation should satisfy the following conditions
\begin{equation}
\label{eq:dv_cond}
    \langle\delta V_i,\phi_j\rangle_{\mathcal M}
    +
    \langle\phi_i,\delta V_j\rangle_{\mathcal M}
    =0,
    \qquad 1\le i,j\le k.
\end{equation}
Using the orthogonal decomposition $T_{[x]^*}\mathcal M=W^-\oplus W^+$, write
\begin{equation}
    \delta V_i
    =
    \sum_{j=1}^k c_{ij}\phi_j+w_i,
    \qquad
    w_i\in W^+.
\end{equation}
Applying Equation \eqref{eq:dv_cond}, we have $c_{ii}=0$ and $c_{ij}+c_{ji}=0$. Thus the coefficients $c_{ij}$ form a skew-symmetric matrix. For $p<q$, set $a_{pq}:=c_{pq}$. Then every tangent variation can be written uniquely as
\begin{equation}\label{eq:frame_tangent_decomposition}
    \delta V_i
    =
    -\sum_{1 \leq p<i}a_{pi}\phi_p
    +
    \sum_{i<q\le k}a_{iq}\phi_q
    +
    w_i,
    \qquad
    w_i\in W^+.
\end{equation}

We now linearize the frame dynamics in Equation \eqref{eq:new_intrinsic_hisd}. With $[x]=[x]^*$ fixed, the $i$th frame dynamic is
\begin{equation}
\label{eq:i-frame}
    \frac{\mathrm dV_i}{\mathrm dt} = F_i([x]^*,\mathbf V)
    =
    -\mathcal P_i(\mathbf V)LV_i.
\end{equation}
The linearized dynamic of Equation \eqref{eq:i-frame} at $\mathbf V^*$ is
\begin{equation}\label{eq:frame_full_linearization}
    \frac{\mathrm d \delta V_i}{\mathrm dt}
    =
    \sum_{p<i}(\lambda_i-\lambda_p)a_{pi}\phi_p
    +
    \sum_{i<q\le k}(\lambda_i-\lambda_q)a_{iq}\phi_q
    -
    (L_+-\lambda_iI)w_i,
\end{equation}
where $L_+:=L|_{W^+}$. This decomposition separates the linearized dynamics into variations within $W^-$ and variations in $W^+$.

For each $p<q$, the coefficient $a_{pq}$ appears in the $p$th and $q$th
components as
$a_{pq}\phi_q$ and $-a_{pq}\phi_p$, respectively.
Equation \eqref{eq:frame_full_linearization} therefore gives
\begin{equation}
    \frac{\mathrm da_{pq}}{\mathrm dt}
    =
    (\lambda_p-\lambda_q)a_{pq}.
\end{equation}
Since $\lambda_p<\lambda_q$, $a_{pq}$ decays
exponentially.
The components in $W^+$ satisfy
\begin{equation}
    \frac{\mathrm dw_i}{\mathrm dt}
    =
    -(L_+-\lambda_iI)w_i.
\end{equation}
Since $[x]^*$ is nondegenerate of index $k$, $W^+$ is the positive spectral
subspace of $L$. By Assumption \ref{hyp:transversal_morse},
$\sigma(L_+)\subset[\gamma_+,\infty)$, and hence
\begin{equation}
    \sigma\bigl(-(L_+-\lambda_iI)\bigr)
    \subset
    (-\infty,-(\gamma_+-\lambda_i)].
\end{equation}
Thus every $\left\|w_i\right\|_{V^*}$ also decays exponentially.

Consequently, there exists $\alpha_V>0$ such that every solution of the linearized frame equation satisfies
\begin{equation}
\left\| \delta\mathbf V(t)\right\|_{V^*}
\le
\exp({-\alpha_Vt})
\left\|\delta\mathbf V(0)\right\|_{V^*}.
\end{equation}
Therefore, $\mathbf V^*$ is a linearly exponentially stable equilibrium of the frame dynamics in $\operatorname{St}_k(T_{[x]^*}\mathcal M)$.
\end{proof}

\begin{theorem}
\label{thm:new_intrinsic_stability}
Under Assumptions \ref{hyp:regular_quotient}, \ref{hyp:transversal_morse}, and \ref{hyp:local_semiflow}, let $[x]^*\in\mathcal M$ be a nondegenerate index-$k$ critical point, and assume that the $k$ negative eigenvalues of $L:=\operatorname{Hess}\mathcal E([x]^*)$ are ordered as $\lambda_1<\cdots<\lambda_k<0,$ with the corresponding orthonormal eigenframe $\mathbf V^*=(\phi_1,\dots,\phi_k)$. Then the equilibrium $([x]^*,\mathbf V^*)$ of Equation \eqref{eq:new_intrinsic_hisd} is locally exponentially stable on $\operatorname{St}_k(T\mathcal M)$.
\end{theorem}

\begin{proof}
The tangent variation of Equation \eqref{eq:quotient_hisd_vector_field} can be presented by
\begin{equation}
\label{eq:cartpro}
    (\delta x,\delta\mathbf V)
    \in
    T_{[x]^*}\mathcal M
    \times
    T_{\mathbf V^*}\operatorname{St}_k(T_{[x]^*}\mathcal M), 
\end{equation}
where the norm of $T_{[x]^*}\mathcal M$ is defined as $\left\|\cdot \right\|_{[x]^*}$ and $\left\| (\delta x,\delta\mathbf V)\right\|_* = \sqrt{\left\|\delta x \right\|_{[x]^*}^2 + \left\| \delta\mathbf V\right\|_{V^*}^2}$ is the joint norm on the Cartesian product \eqref{eq:cartpro}. 

In these coordinates, the Fr\'echet derivative of $F$ at the equilibrium
takes the block form
\begin{equation}
    DF([x]^*,\mathbf V^*)
    =
    \left.
    \left(
    \begin{array}{c:c}
        D_x F_x & D_{\mathbf V}F_x \\[2pt]
        \hdashline
        D_x F_1 & D_{\mathbf V}F_1 \\
        \vdots  & \vdots \\
        D_x F_k & D_{\mathbf V}F_k
    \end{array}
    \right)
    \right|_{([x]^*,\mathbf V^*)}
    =
    \begin{pmatrix}
        J_{xx} & J_{x\mathbf V}\\
        J_{\mathbf Vx} & J_{\mathbf V\mathbf V}
    \end{pmatrix}.
\end{equation}

We first consider the coupling from the frame variables to the position
variable. 
For any frame variation $\delta\mathbf V$ we have
\begin{equation}
    J_{x\mathbf V}\delta\mathbf V
    =
    -D_{\mathbf V}\mathcal A(\mathbf V^*)[\delta\mathbf V]\,
    \operatorname{grad}\mathcal E([x]^*).
\end{equation}
Since $[x]^*$ is a critical point of $\mathcal E$, we have
$\operatorname{grad}\mathcal E([x]^*)=0$ and therefore $J_{x\mathbf V}=0$.
Hence, the full linearized system is block lower triangular.

We next examine the base block. Differentiating $F_x$ with respect to the base variable at fixed
$\mathbf V^*$ gives $J_{xx}\delta x
    =
    -\mathcal A(\mathbf V^*)L\,\delta x$. Since $L\phi_i=\lambda_i\phi_i$, we have $\mathcal A(\mathbf V^*)\phi_i=-\phi_i$.
Thus,
\begin{equation}
    J_{xx}\phi_i
    =
    -\mathcal A(\mathbf V^*)L\phi_i
    =
    \lambda_i\phi_i,
    \qquad
    i=1,\dots,k.
\end{equation}
Since $\lambda_i<0$, $\phi_i$ mode of $\delta x$  decays exponentially for $i=1,\dots, k$.
By Assumption \ref{hyp:transversal_morse} and Lemma \ref{lem:frame_linearization}, $W^+$ is the positive spectral subspace of $L$. On $W^+$,
$\mathcal A(\mathbf V^*)$ is the identity, and therefore $J_{xx}|_{W^+}=-L_+$.
The spectral-gap condition gives
$\sigma(L_+)\subset[\gamma_+,\infty)$, while the negative eigenvalues satisfy
$\lambda_i\le-\gamma_-$. Hence the spectrum of $J_{xx}$ lies in $(-\infty,-\alpha_x]$, where $\alpha_x:=\min\{\gamma_-,\gamma_+\}>0$,
and the spectral theorem yields
\begin{equation}
    \left\|\delta x(t)\right\|_{[x]^*}
    \le
    \exp({-\alpha_x t})\left\|\delta x(0)\right\|_{[x]^*}.
\end{equation}

We now turn to the frame block. 
By definition,
$J_{\mathbf V\mathbf V}$ is obtained by linearizing the direction equations
with respect to $\mathbf V$ while keeping the base point fixed at $[x]^*$.
It is therefore precisely the fixed-base frame linearization considered in
Lemma \ref{lem:frame_linearization}.
Hence, there exist constants
$C_V\ge1$ and $\alpha_V>0$ such that
\begin{equation}
    \left\|\exp({tJ_{\mathbf V\mathbf V}})\right\|
    \le
    C_V \exp({-\alpha_Vt}).
\end{equation}
Consequently, the full linearized system can be written as
\begin{equation}
    \left\{
    \begin{aligned}
        \frac{\mathrm d\delta x}{\mathrm dt}
        &=
        J_{xx}\delta x,\\
        \frac{\mathrm d\delta{\mathbf V}}{\mathrm dt}
        &=
        J_{\mathbf Vx}\delta x
        +
        J_{\mathbf V\mathbf V}\delta\mathbf V.
    \end{aligned}
    \right.
\end{equation}
The first equation gives $\|\delta x(t)\|_{[x]^*} \le \exp({-\alpha_x t})\|\delta x(0)\|_{[x]^*}$. For the frame variables, variation of constants yields
\begin{equation}
    \delta\mathbf V(t)
    =
    \exp({tJ_{\mathbf V\mathbf V}})\delta\mathbf V(0)
    +
    \int_0^t
    \exp({(t-s)J_{\mathbf V\mathbf V}})
    J_{\mathbf Vx}\delta x(s)\,\mathrm ds.
\end{equation}
Since the vector field is $C^1$ in the chosen chart by Assumption \ref{hyp:local_semiflow},
$J_{\mathbf Vx}$ is bounded. Denoting
$M:=\left\|J_{\mathbf Vx}\right\|$, we obtain
\begin{equation}
\begin{aligned}
    \left\|\delta\mathbf V(t)\right\|_{V^*}
    &\le
    C_Ve^{-\alpha_Vt}\left\|\delta\mathbf V(0)\right\|_{V^*}\\
    &\quad+
    C_VM
    \int_0^t
    \exp({-\alpha_V(t-s)})
    \exp({-\alpha_xs})\,\mathrm ds\,
    \left\|\delta x(0)\right\|_{[x]^*}.
\end{aligned}
\end{equation}
The convolution of the two exponentially decaying factors is again
exponentially decaying. More precisely, for any
$0<\alpha<\min\{\alpha_x,\alpha_V\}$, there exists $C_\alpha>0$ such that
\begin{equation}
    \int_0^t
    \exp({-\alpha_V(t-s)})
    \exp({-\alpha_xs})\,\mathrm ds
    \le
    C_\alpha \exp({-\alpha t}).
\end{equation}
Therefore, for suitable constants $C>0$ and $\alpha>0$,
\begin{equation}
    \left\|(\delta x(t),\delta\mathbf V(t))\right\|_*
    \le
    C \exp({-\alpha t})
    \left\|(\delta x(0),\delta\mathbf V(0))\right\|_*.
\end{equation}
Thus, the linearized evolution of \eqref{eq:new_intrinsic_hisd} at
$([x]^*,\mathbf V^*)$ is exponentially stable.

Finally, Assumption \ref{hyp:local_semiflow} allows us to apply the standard
principle of linearized stability for $C^1$ semiflows on Banach spaces.
Hence $([x]^*,\mathbf V^*)$ is locally exponentially stable on
$\operatorname{St}_k(T\mathcal M)$
\cite{henry2006geometric,pazy2012semigroups, smoller2012shock}.
\end{proof}

We next consider critical points whose index differs from $k$.

\begin{theorem}
\label{thm:index_mismatch_instability}
Under Assumptions \ref{hyp:regular_quotient},
\ref{hyp:transversal_morse}, and \ref{hyp:local_semiflow},
let $[x]^*\in\mathcal M$ be a nondegenerate index-$\ell$
critical point, where $k,\ell\in\mathbb N_0$ and $\ell\neq k$.
If $([x]^*,\mathbf V^*)$ is an equilibrium of
Equation \eqref{eq:new_intrinsic_hisd}, where
$\mathbf V^*=(\psi_1,\ldots,\psi_k)$, then it is linearly
unstable and is not locally asymptotically stable on
$\operatorname{St}_k(T\mathcal M)$.
\end{theorem}

\begin{proof}
Let
\begin{equation}
    L:=\operatorname{Hess}\mathcal E([x]^*),
    \qquad
    S:=\operatorname{span}\{\psi_1,\ldots,\psi_k\},
    \qquad
    \mu_i:=\langle L\psi_i,\psi_i\rangle_{\mathcal M}.
\end{equation}
Let $W^-$ be the $\ell$--dimensional subspace spanned by the negative eigenvectors of $L$.
When $k=0$, we set $S=\{0\}$ and omit the frame variables.

We use the block notation and the joint norm $\left\|\cdot\right\|_*$
introduced in the proof of
Theorem \ref{thm:new_intrinsic_stability}.
Choose local coordinates $(\xi,\eta)$ centered at
$([x]^*,\mathbf V^*)$, where $\xi$ and $\eta$ are the position
and frame coordinates, respectively.

We first consider the coupling from the frame variables to
the position variable.
Since $\operatorname{grad}\mathcal E([x]^*)=0$, the same
calculation as in Theorem \ref{thm:new_intrinsic_stability} gives
\begin{equation}
    J_{x\mathbf V}=0,
    \qquad
    J_{xx}=-\mathcal A(\mathbf V^*)L.
\end{equation}
Hence, the full linearized system is
\begin{equation}\label{eq:mismatch_linearized_system}
    \left\{
    \begin{aligned}
        \frac{\mathrm d\delta x}{\mathrm dt}
        &=J_{xx}\delta x,\\
        \frac{\mathrm d\delta{\mathbf V}}{\mathrm dt}
        &=J_{\mathbf Vx}\delta x
        +J_{\mathbf V\mathbf V}\delta\mathbf V.
    \end{aligned}
    \right.
\end{equation}

We next examine the base block.
The equilibrium conditions for the frame variables give
\begin{equation}
    L\psi_i
    =
    \mu_i\psi_i
    +
    2\sum_{j=1}^{i-1}
    \langle L\psi_i,\psi_j\rangle_{\mathcal M}\psi_j,
    \qquad i=1,\ldots,k.
\end{equation}
Taking the inner product with $\psi_j$ for $j<i$ yields
$\langle L\psi_i,\psi_j\rangle_{\mathcal M}=0$.
Therefore, we have $L\psi_i=\mu_i\psi_i, \forall i=1,\ldots,k$, i.e., $(\mu_i,\psi_i)$ is an eigenpair of $L$ and thus $\mu_i \neq 0$. Since $L$ is self-adjoint, both $S$ and $S^\perp$ are invariant
under $L$.
By the definition of $\mathcal A(\mathbf V^*)$, we have $J_{xx}|_S=L|_S$ and $J_{xx}|_{S^\perp}=-L|_{S^\perp}$.
Thus $J_{xx}$ is self-adjoint.

If $\ell>k$, then
$\dim(W^-\cap S^\perp)\ge\ell-k>0$.
This intersection is finite-dimensional and invariant under $L$,
so it contains an eigenvector $u$ with
$Lu=\lambda u$ and $\lambda<0$.
Consequently, $J_{xx}u=-\lambda u$.

If $\ell<k$, at least one $\mu_i$ is positive.
Otherwise, all $k$ orthonormal vectors $\psi_i$ would belong
to $W^-$, contradicting $\dim W^-=\ell<k$.
For such an index $i$, we have
$J_{xx}\psi_i=\mu_i\psi_i$.

In both cases, there exist $r>0$ and a unit vector $u$ such
that $J_{xx}u=ru$.
We now consider the solution of
Equation \eqref{eq:mismatch_linearized_system} with
$\delta x(0)=u$ and $\delta\mathbf V(0)=0$.
By the self-adjointness of $J_{xx}$,
\begin{equation}
\label{eq:positive_eigen}
    \frac{\mathrm d}{\mathrm dt}
    \langle\delta x(t),u\rangle_{\mathcal M}
    =
    \langle J_{xx}\delta x(t),u\rangle_{\mathcal M}
    =
    r\langle\delta x(t),u\rangle_{\mathcal M}.
\end{equation}
Since $\langle\delta x(0),u\rangle_{\mathcal M}=1$,
integration gives
\begin{equation}
    \left\|(\delta x(t),\delta\mathbf V(t))\right\|_*
    \ge
    \left\|\delta x(t)\right\|_{[x]^*}
    \ge
    \left|\langle\delta x(t),u\rangle_{\mathcal M}\right|
    =
    e^{rt}.
\end{equation}
Therefore the full linearized system is unstable.

We finally establish nonlinear instability by applying the
principle of linearized instability
\cite{daleckii2002stability,ruelle2014elements}.
Let $J:=DF([x]^*,\mathbf V^*)$.
Under Assumption \ref{hyp:local_semiflow}, the vector field
is $C^1$ in the chosen local coordinates, and $J$ is bounded.
The principle states that the equilibrium is unstable if
there exists $a\ge0$ such that
\begin{equation}\label{eq:mismatch_spectral_separation}
\begin{aligned}
    \sigma(J)\cap
    \{\lambda\in\mathbb C:\operatorname{Re}\lambda=a\}
    =\varnothing,\qquad
    \sigma(J)\cap
    \{\lambda\in\mathbb C:\operatorname{Re}\lambda>a\}
    \neq\varnothing.
\end{aligned}
\end{equation}
We now verify this spectral condition.

For $k\ge1$, the same calculation as in
Lemma \ref{lem:frame_linearization} shows that
$J_{\mathbf V\mathbf V}$ is self-adjoint and
\begin{equation}
    \sigma(J_{\mathbf V\mathbf V})
    =
    \{\mu_p-\mu_q:1\le p<q\le k\}
    \cup
    \{\mu_i-\lambda:1\leq i \leq k, \lambda\in\sigma(L|_{S^\perp})\}.
\end{equation}
Since $\sigma(J_{xx}) = \{\mu_1,\dots,\mu_k\} \cup \sigma(-L|_{S^\perp})$ and both diagonal blocks are self-adjoint, the block
lower-triangular form gives
\begin{equation}
    \sigma(J)
    =
    \sigma(J_{xx})\cup\sigma(J_{\mathbf V\mathbf V})
    \subset\mathbb R.
\end{equation}

Since $S+W^-$ is finite-dimensional and $L\ge\gamma_+I$ on $(S+W^-)^\perp$, $\sigma(-L|_{S^\perp\cap W^-})$ is a finite set and $\sigma(-L|_{(S+W^-)^\perp})
\subset(-\infty,-\gamma_+]$. By the decomposition $\sigma(-L|_{S^\perp}) = \sigma(-L|_{S^\perp\cap W^-})\cup \sigma(-L|_{(S+W^-)^\perp})$, we take a finite subset of $\sigma(J)$ as
\begin{equation}
\begin{aligned}
    \Sigma_f = & \{\mu_1,\dots,\mu_k\} \cup \sigma(-L|_{S^\perp\cap W^-}) \\
    &\cup \{\mu_p-\mu_q:1\le p<q\le k\}
    \cup\{\mu_i-\lambda:1\leq i \leq k, \lambda\in\sigma(L|_{S^\perp\cap W^-})\}.
    \end{aligned}
\end{equation}
Denote $\mu_{\max} = \max\{\mu_1,\dots,\mu_k\} \neq 0$, and the rest part of $\sigma(J)$ can be bounded as
\begin{equation}
\begin{aligned}    
\sigma(J)\setminus\Sigma_f &\subset \sigma(-L|_{(S+W^-)^\perp})\cup
    \{\mu_i-\lambda:1\leq i \leq k, \lambda\in\sigma(L|_{(S+W^-)^\perp})\} \\
    &\subset (-\infty, \max\{0,\mu_{\max}\}-\gamma_+].
\end{aligned}
\end{equation}
If $\mu_{\max} > 0$, then we can take $a \in (\mu_{\max}-\gamma_+, \mu_{\max})\setminus \Sigma_f$ to satisfy Equation \eqref{eq:mismatch_spectral_separation}. 
Otherwise, $\mu_{\max} < 0$ indicates that $r \in \Sigma_f$ in Equation \eqref{eq:positive_eigen} is an isolated positive eigenvalue of $J$.
Therefore, the nonlinear instability is obtained.

For $k=0$, we have $J=-L$.
Since $\ell>0$ and $W^-$ is finite-dimensional, the positive
spectrum of $J$ is a nonempty finite set, while its remaining
spectrum lies in $(-\infty,-\gamma_+]$.
Hence the largest spectral value is again positive and isolated.

In both cases, choosing $a>0$ sufficiently close to the
largest spectral value from below gives
Equation \eqref{eq:mismatch_spectral_separation}.
The principle of linearized instability therefore applies,
and $([x]^*,\mathbf V^*)$ is unstable.
In particular, it is not locally asymptotically stable on
$\operatorname{St}_k(T\mathcal M)$.
\end{proof}

\subsection{SR-HiSD in the ambient space} 
We next express the quotient HiSD \eqref{eq:new_intrinsic_hisd} in the ambient Hilbert space $H$ under Assumption \ref{hyp:regular_quotient}. 
For a quotient trajectory $([x](t),\mathbf V(t))$, 
we aim to obtain a horizontal representative $x(t)$ of $[x](t)$ and a horizontal lift $\mathbf v(t)=(v_1(t),\ldots,v_k(t))$ of $\mathbf V(t)$, satisfying
$\dot x(t)\in\mathcal H_{x(t)}$.

We first derive the $x$-dynamics.
By Equation \eqref{eq:reduced_gradient}, the horizontal lift
of $\operatorname{grad}\mathcal E([x])$ is $\nabla E(x)$.
Since $D\pi_x:\mathcal H_x\to T_{[x]}\mathcal M$ is an isometry,
the first equation in \eqref{eq:new_intrinsic_hisd} becomes
\begin{equation}\label{eq:sr_hisd_lifted_position}
    \frac{\mathrm dx}{\mathrm dt}
    =
    -\bigl(
        \mathcal A(\mathbf V)\operatorname{grad}\mathcal E([x])
    \bigr)^\sharp
    =
    -\left(I-2\sum_{j=1}^k v_jv_j^\top\right)\nabla E(x).
\end{equation}

We next obtain the $v_i$-dynamics.
The Levi--Civita connection on the quotient satisfies $
    \left(\frac{\mathrm DV_i}{\mathrm dt}\right)^\sharp
    =
    P_x\frac{\mathrm dv_i}{\mathrm dt}$ along a horizontal representative
\cite{o1966fundamental}.
On the other hand, differentiating the constraint $P_xv_i=v_i$
along $x(t)$ gives
\begin{equation}\label{eq:sr_hisd_horizontal_derivative}
    DP_x\left(\frac{\mathrm dx}{\mathrm dt}\right)v_i+P_x\frac{\mathrm dv_i}{\mathrm dt}=\frac{\mathrm dv_i}{\mathrm dt}.
\end{equation}
By Equation \eqref{eq:reduced_hessian}, we have 
\begin{equation}
    \bigl(\operatorname{Hess}\mathcal E([x])V_i\bigr)^\sharp
    =
    P_x\nabla^2E(x)P_xv_i.
\end{equation}
Therefore, the horizontal lift of the direction equation in
\eqref{eq:new_intrinsic_hisd} is
\begin{equation}\label{eq:sr_hisd_lifted_direction}
    P_x\frac{\mathrm dv_i}{\mathrm dt}
    =
    -\Bigl(
        I-v_iv_i^\top-2\sum_{j=1}^{i-1}v_jv_j^\top
    \Bigr)
    P_x\nabla^2E(x)P_xv_i,
    \qquad i=1,\ldots,k.
\end{equation}
Combining Equations \eqref{eq:sr_hisd_lifted_position},
\eqref{eq:sr_hisd_horizontal_derivative}, and
\eqref{eq:sr_hisd_lifted_direction}, we obtain SR-HiSD for locating an index-$k$ saddle point ($k$-SR-HiSD):
\begin{equation}\label{eq:new_sr_hisd}
\left\{
\begin{aligned}
    \frac{\mathrm dx}{\mathrm dt}
    &=
    -\Bigl(I-2\sum_{j=1}^k v_jv_j^\top\Bigr)\nabla E(x),
    \\
    \frac{\mathrm dv_i}{\mathrm dt}
    &=
    -\Bigl(
        I-v_iv_i^\top-2\sum_{j=1}^{i-1}v_jv_j^\top
    \Bigr)
    P_x\nabla^2E(x)P_xv_i
    +DP_x\left(\frac{\mathrm dx}{\mathrm dt}\right)v_i,
    \quad 1\leq i \leq k.
\end{aligned}
\right.
\end{equation}
Here the Hessian acts on the horizontal space through
$P_x\nabla^2E(x)P_x$, while the term $DP_x\left(\frac{\mathrm dx}{\mathrm dt}\right)v_i$
accounts for a parallel transport of the horizontal space along $x(t)$.
This term is required to maintain the constraint $P_xv_i=v_i$.

\begin{theorem}
\label{thm:new_ambient_equivalence}
Under Assumption \ref{hyp:regular_quotient}, let
$([x](t),\mathbf V(t))_{t\in [0,\infty)}$ be a $C^1$ solution of
the quotient HiSD \eqref{eq:new_intrinsic_hisd}.
For every $x_0\in \pi^{-1}([x](0))$, there exists a unique
$C^1$ horizontal representative $x(t)$ with $x(0)=x_0$.
Along $x(t)$, the frame $\mathbf V(t)$ has a unique $C^1$
horizontal lift $\mathbf v(t)$.
$(x(t),\mathbf v(t))$ satisfies SR-HiSD \eqref{eq:new_sr_hisd},
with $P_xv_i=v_i$ and
$\langle v_i,v_j\rangle_H=\delta_{ij}$ for all $t\in [0,\infty)$.
\end{theorem}

\begin{proof}
We recall the horizontal lifting result \cite{kriegl1997convenient}.
Let $\pi:\Omega\to\mathcal M$ be a smooth principal $G$-bundle,
where $\mathcal M$ is a Hilbert manifold and $G$ is a
finite-dimensional Lie group.
Assume that the horizontal spaces $\mathcal H_x$ define
a $C^1$ principal connection. For a $C^1$ curve $c:I\to\mathcal M$ with $0\in I$ and
a prescribed point $x_0\in\pi^{-1}(c(0))$, there exists
a unique $C^1$ curve $x:I\to\Omega$ satisfying
$\pi(x(t))=c(t)$, $x(0)=x_0$, and
$\dot x(t)\in\mathcal H_{x(t)}$ for all $t\in I$.
This curve is the horizontal lift of $c(t)$ through $x_0$. The result for smooth connections and curves is given in
\cite{kriegl1997convenient}.
For finite-dimensional $G$, the same proof applies to
a $C^1$ connection and a $C^1$ curve, yielding the
$C^1$ lift stated above.

We now verify that the horizontal lifting result applies to our quotient settings.
Under Assumption \ref{hyp:regular_quotient},
$\pi:\Omega\to\mathcal M$ is a smooth principal $G$-bundle,
and the horizontal spaces $\mathcal H_x$ depend $C^1$ on $x$. Moreover, the $G$-action preserves the horizontal spaces. 
Hence $\mathcal H_x$ defines a $C^1$ principal connection
\cite{kriegl1997convenient}. For the given $C^1$ curve $[x](t)$ and any
$x_0\in\pi^{-1}([x](0))$, the horizontal lifting result yields
a unique $C^1$ curve $x:[0,\infty)\to\Omega$ such that
$\pi(x(t))=[x](t)$, $x(0)=x_0$, and
$\dot x(t)\in\mathcal H_{x(t)}$.

For the given horizontal representative $x(t)$, each $V_i(t)$ has the unique horizontal lift
$v_i(t):=V_i(t)^{\sharp} = (D\pi_{x(t)}|_{\mathcal H_{x(t)}})^{-1}V_i(t)$,
$i=1,\ldots,k$. By Assumption \ref{hyp:regular_quotient},
$D\pi_x|_{\mathcal H_x}$ is an isometric isomorphism.
Since $\pi$ is smooth and $x\mapsto P_x$ is $C^1$,
these isomorphisms depend $C^1$ on $x$ and so do $(D\pi_{x}|_{\mathcal H_{x}})^{-1}$. The inverses are bounded and linear, so the
horizontal lift $(D\pi_{x}|_{\mathcal H_{x}})^{-1}V$ depends $C^1$ on $(x,V)$.
Since $x(t)$ and $V_i(t)$ are $C^1$, the chain rule gives
$v_i\in C^1([0,\infty),H)$. 

Hence the identities
\eqref{eq:sr_hisd_lifted_position},
\eqref{eq:sr_hisd_horizontal_derivative}, and
\eqref{eq:sr_hisd_lifted_direction}, derived above for
$C^1$ horizontal lifts, apply to
$(x(t),\mathbf v(t))$.
Combining these identities gives the SR-HiSD system
\eqref{eq:new_sr_hisd}. 

By construction, $v_i(t)\in\mathcal H_{x(t)}$,
so $P_xv_i=v_i$ for all $t\in[0,\infty)$.
Since $D\pi_x|_{\mathcal H_x}$ is an isometry
and $\mathbf V(t)$ is orthonormal, we have $\langle v_i,v_j\rangle_H = \langle V_i,V_j\rangle_{\mathcal M} = \delta_{ij}$, $1\le i,j\le k$.
Thus the lifted solution satisfies both the
horizontal and orthonormality constraints. The preservation of constraints can also be verified directly by checking that the SR-HiSD vector field is tangent to the corresponding constraint manifold.
\end{proof}

\begin{remark}
\label{rem:quotient_distance}
Under Assumption \ref{hyp:regular_quotient}, the vertical spaces
$\mathcal V_x=T_x(G\cdot x)$ form a smooth finite-dimensional subbundle, and
the corresponding horizontal projector $P_x$ inherits the regularity of the
group action. With sufficient regularity of the group action and $E$, the maps $DP_x$ and
$P_x\nabla^2E(x)P_x$ have the local regularity required above.
\end{remark}

To formulate convergence in the ambient space, let
$\mathcal U_{([x^*],\mathbf V^*)}$ be a local stable neighborhood of the
equilibrium $([x^*],\mathbf V^*)$ in the Stiefel bundle
$\operatorname{St}_k(T\mathcal M)$. We define the corresponding lifted basin
by requiring both the base point and the horizontal frame to map into this
neighborhood:
\begin{equation}
  \mathcal B(x^*,\mathbf v^*)
  :=
  \left\{
    (x,v_1,\dots,v_k):
    \bigl(\pi(x),D\pi_xv_1,\dots,D\pi_xv_k\bigr)
    \in \mathcal U_{([x^*],\mathbf V^*)}
  \right\}.
\end{equation}

Moreover, since the action of $G$ on the Hilbert space $H$ is isometric, the
quotient distance admits the orbit-distance bound
\begin{equation}\label{eq:quotient_orbit_distance}
    d_{\mathcal M}(\pi(x),[y])
        \geq
    \inf_{g\in G}\left\|x-g\cdot y\right\|_H
    =
    \operatorname{dist}_H(x,G\cdot y).
\end{equation}

\begin{corollary}
\label{cor:new_orbital_convergence}
Under the assumptions of Theorem \ref{thm:new_intrinsic_stability}, if the horizontal lift of the initial
state and frame belongs to $\mathcal B(x^*,\mathbf v^*)$, then the corresponding ambient trajectory of
Equation \eqref{eq:new_sr_hisd} satisfies
\begin{equation}
  \operatorname{dist}_H\bigl(x(t),G\cdot x^*\bigr)\to 0
  \qquad \text{as } t\to\infty.
\end{equation}
\end{corollary}

\begin{proof}
The quotient trajectory converges to $[x^*]$ by Theorem
\ref{thm:new_intrinsic_stability}. Using Equation \eqref{eq:quotient_orbit_distance}, we obtain
\begin{equation}
  \operatorname{dist}_H(x(t),G\cdot x^*)
  \leq d_{\mathcal M}([x](t),[x^*])\longrightarrow0.
\end{equation}
\end{proof}

\begin{theorem}[Pointwise Convergence in the Ambient Space]
\label{thm:new_pointwise_convergence}
Under Assumptions \ref{hyp:regular_quotient}--\ref{hyp:local_semiflow}, $(x(t),\mathbf{v}(t))$ is the solution of SR-HiSD \eqref{eq:new_sr_hisd} with an initial condition $(x(0),\mathbf{v}(0))\in \mathcal B(x^*,\mathbf v^*)$. 
Then there exists a point $x_\infty\in G\cdot x^*$ such that
\begin{equation}
  \lim_{t\to\infty}\left\|x(t)-x_\infty\right\|_H = 0.
\end{equation}
In particular, the SR-HiSD position trajectory converges pointwise to an index-$k$ critical point on the critical orbit.
\end{theorem}

\begin{proof}
By Theorem \ref{thm:new_intrinsic_stability}, the trajectory $([x](t),\mathbf V(t))$ converges exponentially to $([x^*],\mathbf V^*)$ on the quotient orthonormal $k$-frame bundle. This nonlinear exponential stability follows from the exponential stability of the linearized evolution and the standard principle
of linearized stability
\cite{henry2006geometric,smoller2012shock}.
Since the quotient HiSD vector field is $C^1$ and vanishes at the equilibrium,
its value along an exponentially convergent trajectory also decays
exponentially. In particular, for the position component there exist constants
$C,\beta>0$ such that
\begin{equation}
  \left\|\tfrac{\mathrm d}{\mathrm dt}[x](t)\right\|_{[x](t)} \le Ce^{-\beta t}.
\end{equation}

Since $x(t)$ is the horizontal lift of $[x](t)$ and $\pi$ is a Riemannian submersion, horizontal vectors preserve norm under $D\pi_x$. Therefore, we have
\begin{equation}
  \left\|\tfrac{\mathrm d}{\mathrm dt}x(t)\right\|_H = \left\|\tfrac{\mathrm d}{\mathrm dt}[x](t)\right\|_{[x](t)} \le Ce^{-\beta t}.
\end{equation}
It follows that the arc length of $x(t)$ is finite:
\begin{equation}
  \int_0^\infty \left\|\tfrac{\mathrm d}{\mathrm dt}x(t)\right\|_H\,dt < \infty.
\end{equation}
Consequently, for any $t \geq 0$, $\Delta t >0$, we have 
\begin{equation}
  \left\|x(t+\Delta t)-x(t)\right\|_H
  \le \left\|\int_{t}^{t+\Delta t}\tfrac{\mathrm d}{\mathrm dt}x(t)\,dt\right\|_H \le \int_{t}^{\infty}Ce^{-\beta s}\,ds
  \le \frac{C}{\beta}e^{-\beta t},
\end{equation}
Hence, $\{x(t)\}_{t\geq 0}$ satisfies the Cauchy criterion as $t \to +\infty$, and $\displaystyle\lim_{t\to +\infty}x(t)=x_{\infty}$ exists.

By Corollary \ref{cor:new_orbital_convergence}, the distance from $x(t)$ to the orbit $G\cdot x^*$ tends to zero. Passing to the limit gives $\operatorname{dist}_H(x_\infty,G\cdot x^*)=0$, which implies that $x_\infty\in G\cdot x^*$. Thus the ambient trajectory converges pointwise to a single point on the critical orbit, which is an index-$k$ critical point.
\end{proof}

\section{Algorithm}
\label{sec:algorithm_q}
\subsection{Numerical implementation of the horizontal projector}
The horizontal projector $P_x$ to the horizontal space $\mathcal H_x$ can be computed using the
infinitesimal generators of the group action.
Let $\mathfrak g=T_eG$ be the Lie algebra of $G$,
with a fixed basis $\{\xi_1,\ldots,\xi_m\}$,
where $m=\dim G$.
Under Assumption \ref{hyp:regular_quotient}, for each $x\in\Omega$, vectors
\begin{equation}
    z_p(x)
    :=
    \left.\frac{\mathrm d}{\mathrm dt}\right|_{t=0}
    \bigl(\exp(t\xi_p)\cdot x\bigr),
    \qquad p=1,\ldots,m.
\end{equation}
form a basis of the vertical space $\mathcal V_x$.

In the numerical implementation, we apply Gram--Schmidt orthonormalization to
$z_1(x),\ldots,z_m(x)$ to obtain an orthonormal basis of the tangent space.
In practice, the dynamic may approach states with additional continuous symmetry in $G$, making some orbit tangent vectors nearly linearly dependent.
Therefore, we use a tolerance $\epsilon >0$ to avoid numerical instability caused by small residuals. 
At step $p$, we subtract from $z_p(x)$ its
projections onto the vectors retained in the preceding steps:
\begin{equation}
    r_p
    :=
    z_p(x)
    -
    \sum_{\substack{1\le j<p}}
    \langle z_p(x),q_j(x)\rangle_H q_j(x).
\end{equation}
If $\|r_p\|_H>\epsilon$, we set
$q_p(x):=r_p/\|r_p\|_H$ and add it to the basis.
Otherwise, we proceed to the next step without adding a vector. 

After processing all $m$ vectors, let $s\le m$ be the
number of retained vectors.
We then compute the horizontal projector
$P_x:H\to\mathcal H_x$ by
\begin{equation}\label{eq:numerical_horizontal_projection}
    P_xu
    :=
    u-\sum_{j=1}^{s}
    \langle u,q_j(x)\rangle_H q_j(x),
    \qquad u\in H.
\end{equation}
In practice, the horizontal projector can be updated every several iterations of the position variable $x$ to reduce the computational cost.

\subsection{Numerical algorithms of SR-HiSD}

The numerical scheme for SR-HiSD \eqref{eq:new_sr_hisd} uses
an explicit Euler step for the position and a general eigensolver $\mathtt{EigenSol}_k$
for the directions.
Given an orthonormal horizontal frame at $x^{(n)}$ as $[v_i^{(n)}]_{i=1}^k$, the
iteration reads
\begin{equation}\label{eq:disc_dynamics}
\left\{
\begin{aligned}
    &x^{(n+1)}
    = x^{(n)}-\tau^{(n)}
       \left(
          \nabla E(x^{(n)})
          -2\sum_{j=1}^k
          \big\langle\nabla E(x^{(n)}),v_j^{(n)}\big\rangle_H
          v_j^{(n)}
       \right),\\
    &[v_1^{(n+1)},\ldots,v_k^{(n+1)}]
    = \mathtt{EigenSol}_k\!\left(
        L_{x^{(n+1)}};
        [P_{x^{(n+1)}}v_1^{(n)},\ldots,
         P_{x^{(n+1)}}v_k^{(n)}]
       \right).
\end{aligned}
\right.
\end{equation}
Here $\tau^{(n)}>0$ is the position step size, and
$L_x:=\left.P_x\nabla^2E(x)P_x\right|_{\mathcal H_x}$.
The eigensolver $\mathtt{EigenSol}_k$ uses the projected
previous frame as its initial guess 
and returns vectors which approximate eigenvectors associated with the $k$ smallest eigenvalues of $L_{x^{(n+1)}}$. 
Since the operator $L_{x^{(n+1)}}$ is defined on $\mathcal H_{x^{(n+1)}}$, the eigenproblem is solved only in this space, and the output vectors are orthonormal in $\mathcal H_{x^{(n+1)}}$.  
In numerical experiments, we apply two steps of the Locally Optimal Block Preconditioned Conjugate Gradient (LOBPCG) method at each iteration as the eigensolver \cite{knyazev2001toward}.

The position step size $\tau^{(n)}$ can be fixed or determined
by the Barzilai--Borwein (BB) gradient method \cite{barzilai1988two}.
Let $g^{(n)}$ denote the reflected gradient in the position
update \eqref{eq:disc_dynamics}, and set
$s^{(n)}:=x^{(n)}-x^{(n-1)}$, $y^{(n)}:=g^{(n)}-g^{(n-1)}$.
The BB step size is 
$\tau^{(n)} = \langle s^{(n)},y^{(n)}\rangle_H/\|y^{(n)}\|_H^2$, 
and the practical step size is restricted to
$[\tau_{\min},\tau_{\max}]$ to avoid excessively small or large steps.

\section{Numerical experiments}\label{sec:experiments_q}
\subsection{Morse clusters with rigid-motion symmetry}
We use Morse clusters to compare SR-HiSD with the original HiSD \cite{doye1995effect}.
For an $N$-particle cluster in $\mathbb R^d$ with positions $r=(\mathbf r_1,\ldots,\mathbf r_N)\in(\mathbb R^d)^N$, the total energy is
\begin{equation}
    E(r)
    =
    \sum_{1\le i<j\le N}
    V\left(\|\mathbf r_i-\mathbf r_j\|\right).
\end{equation}
The pairwise interaction $V$ is taken as 
\begin{equation}
    V(R)
    =
    \left(1-\exp({-\rho_0(R-1)})\right)^2-1,
\end{equation}
where $\rho_0>0$ controls the width of the potential well \cite{doye1995effect}. 
This energy is invariant under rigid motions, described by $\mathrm{SE}(d)$. 
At a critical point of $E$ where the action is free,
the orbit tangent space has dimension $m_d={d(d+1)}/{2}$,
with $m_2=3$ and $m_3=6$.
Such symmetry-induced degeneracy also occurs in
Lennard--Jones clusters \cite{wales1997global} and
molecular systems.

To compute the horizontal projector, let $e_\alpha$ be
the standard basis vectors of $\mathbb R^d$ and set
$A_{pq}:=e_qe_p^\top-e_pe_q^\top$.
The translational and rotational directions are
\begin{align}
    z_p^{\rm tr}(r)
    &=(e_p,\ldots,e_p),
    &&p=1,\ldots,d,\label{eqn:translation}\\
    z_{pq}^{\rm rot}(r)
    &=(A_{pq}\mathbf r_1,\ldots,A_{pq}\mathbf r_N),
    &&1\le p<q\le d,\label{eqn:rotation}
\end{align}
which span the vertical space $\mathcal V_r$.
It should be noted that, at a nonstationary point $r$ of $E$, the translational direction $z_p^{\rm tr}(r)$ lies in the kernel of $\nabla^2E(r)$, whereas the rotational direction $z_{pq}^{\rm rot}(r)$ need not.

We compare 1-SR-HiSD with 1-, 2-, 3-, and 4-HiSD and NPSS \cite{cui2025efficient} on two-dimensional $(d=2)$ Morse clusters with $N=5,\ldots,10$ and $\rho_0=6$. 
For each $N$, we generate $500$ initial configurations by perturbing a local minimizer by a distance of $0.10$ along uniformly sampled directions in the horizontal space. 
We obtain $3000$ initial configurations in total, which is applied to each method. 
For NPSS algorithm, the nullspace update threshold is set to $\epsilon_\lambda=0.10$ or $0.01$. 
Each method initializes its ascent directions using the smallest Hessian eigenvectors, with the corresponding horizontal or nullspace constraints applied for SR-HiSD and NPSS. 
All methods use the same BB step-size strategy and LOBPCG-based direction updates, with an initial step size of $\Delta t=0.02$, a maximum of $50000$ iterations, and a stopping tolerance of $10^{-8}$. 
A run is counted as successful if it converges to an index-1 saddle point. 

\begin{table}[htbp]
\centering
\caption{Success rates for reaching index-$1$
saddle points in two-dimensional Morse clusters.}
\label{tab:morse_counterexample_success}
\small
\setlength{\tabcolsep}{6pt}
\begin{tabular}{@{}lcccccc@{}}
\hline
Method & $N=5$ & $N=6$ & $N=7$ & $N=8$ & $N=9$ & $N=10$ \\
\hline
1-HiSD & $2.4\%$ & $0.0\%$ & $0.0\%$ & $0.0\%$ & $0.0\%$ & $0.0\%$ \\
2-HiSD & $3.8\%$ & $3.6\%$ & $3.6\%$ & $0.0\%$ & $0.4\%$ & $0.4\%$ \\
3-HiSD & $4.4\%$ & $3.6\%$ & $3.6\%$ & $0.0\%$ & $0.4\%$ & $0.4\%$ \\
4-HiSD & $92.6\%$ & $99.0\%$ & $69.4\%$ & $85.0\%$ & $45.4\%$ & $90.0\%$ \\
\hline

NPSS, $\epsilon_\lambda=0.10$
    & $89.8\%$ & $91.2\%$ & $82.4\%$ & $99.8\%$ & $98.0\%$ & $98.6\%$ \\

NPSS, $\epsilon_\lambda=0.01$
    & $90.6\%$ & $92.0\%$ & $84.2\%$ & $99.8\%$ & $98.0\%$ & $98.6\%$ \\
\hline
1-SR-HiSD & $\mathbf{100.0\%}$ & $\mathbf{100.0\%}$ & $\mathbf{100.0\%}$ & $\mathbf{100.0\%}$ & $\mathbf{100.0\%}$ & $\mathbf{100.0\%}$ \\
\hline
\end{tabular}
\end{table}

Table~\ref{tab:morse_counterexample_success} reports the success rates. 
Most trials using 1-, 2-, and 3-HiSD converge to local minimizers.  
The success rates of 4-HiSD are much higher, and most unsuccessful searches converge to index-2 saddles. 
For example, $30.0\%$ $(N=7)$ and $52.8\%$ $(N=9)$ of the 4-HiSD runs converge to index-2 saddles. 
The NPSS method substantially improves the success rates, but most of its unsuccessful searches return to the initial local minimizer. 
The critical reason is that the rotational directions $z_{pq}^{\rm rot}(r)$ generally do not lie in the Hessian nullspace when $r$ is away from critical points.
As the iterations progress, the Hessian nullspace cannot capture these directions, which are always orthogonal to the gradient. 
In these unsuccessful searches, the ascent direction is dominated by the rotational directions and becomes nearly orthogonal to the gradient, making the position update close to gradient descent and finally converging to a local minimizer. 
In contrast, SR-HiSD constructs the horizontal projector directly from the translational and rotational directions, so that the ascent directions are always orthogonal to rigid-motion directions. 
As a result, 1-SR-HiSD successfully locates index-1 saddle points in all $3000$ trials, demonstrating its reliability in searching for saddle points of the prescribed Morse index.

\begin{figure}[htbp]
    \centering
    \includegraphics[width=0.95\textwidth]{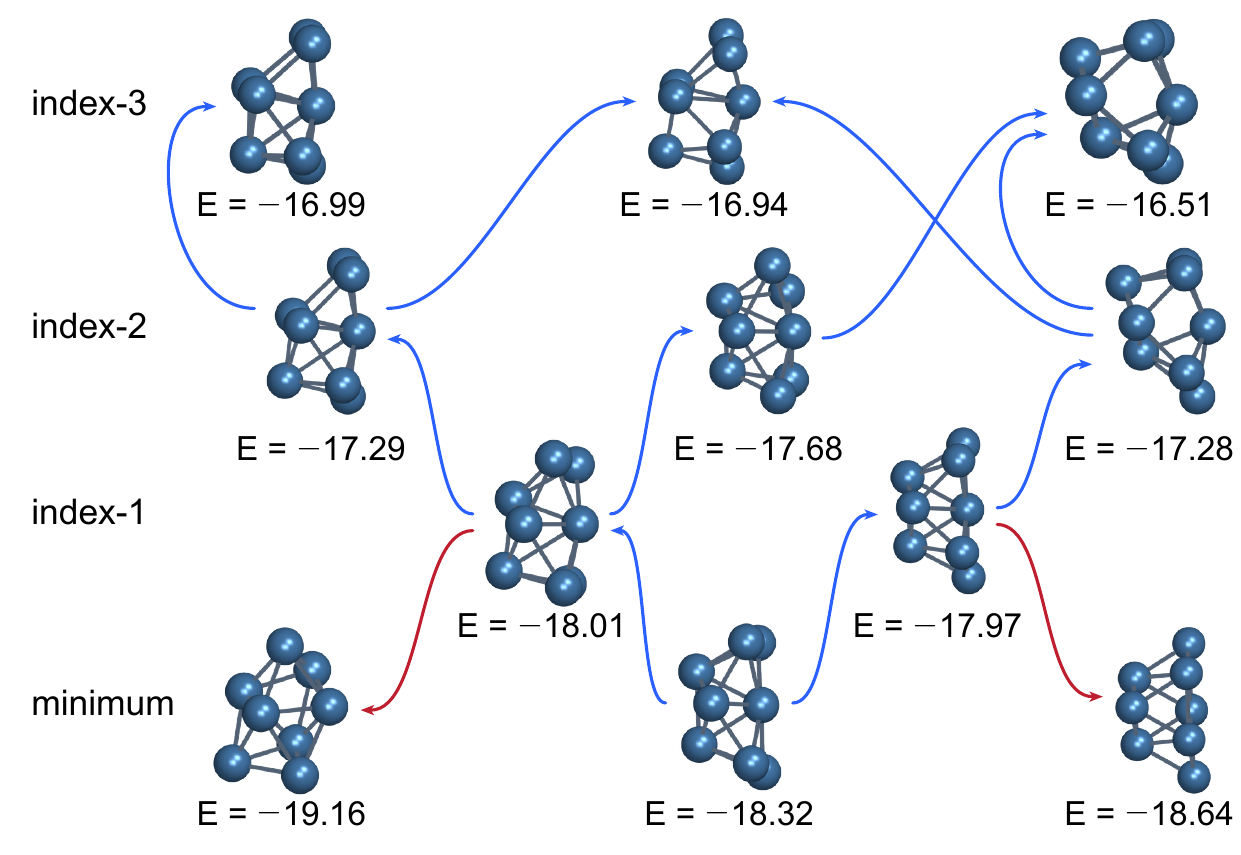}
    \caption{SR-HiSD saddle search for the 3D Morse cluster with $N=8$, starting from multi-start minima and searching upward to higher-index saddles.}
    \label{fig:morse3d_n8_landscape}
\end{figure}

We next use SR-HiSD to construct the solution landscape of a three-dimensional ($d=3$) Morse cluster with $N=8$ and $\rho_0=6$ following the method in \cite{yin2020construction}.
At a regular critical point, the Hessian kernel contains six directions generated by the $\mathrm{SE}(3)$ action,  three translational and three rotational directions.
Starting from a local minimizer with energy $E\approx-18.32$, we perform successive upward searches with SR-HiSD from index-$1$ to index-$3$. 
For each index-1 saddle, we perturb on both sides along the negative horizontal Hessian eigenvector and follow the gradient flow to locate the connected minima.
Figure~\ref{fig:morse3d_n8_landscape} shows the resulting solution landscape, with symmetry-equivalent saddle points represented by a single critical orbit.
The blue arrows indicate upward searches by SR-HiSD, while the red arrows indicate gradient-flow trajectories.
These results demonstrate that SR-HiSD can locate saddles with the prescribed indices and reveal connection relationships among stationary states.

\subsection{Onsager model with rotational symmetry}

We next consider the Onsager model for nematic liquid crystals
on the unit sphere $\mathbb S^2$ \cite{yin2022solution}.
The state variable is an orientational probability density
\begin{equation}
    \rho\in\mathcal P(\mathbb S^2)
    :=
    \left\{
    \rho\ge0:
    \int_{\mathbb S^2}\rho(\mathbf p)\,\mathrm d\mathbf p=1
    \right\}.
\end{equation}
With the temperature coefficient normalized to one,
the free energy is
\begin{equation}
    E_\mu(\rho)
    =
    \int_{\mathbb S^2}
    \rho(\mathbf p)\log\rho(\mathbf p)\,\mathrm d\mathbf p
    +
    \frac12
    \int_{\mathbb S^2}\int_{\mathbb S^2}
    k_\mu(\mathbf p\cdot\mathbf q)
    \rho(\mathbf p)\rho(\mathbf q)\,
    \mathrm d\mathbf p\,\mathrm d\mathbf q,
\end{equation}
where $k_\mu(t)=\mu\sqrt{1-t^2}$ is the Onsager
interaction kernel and $\mu>0$ is the interaction parameter.
The energy is invariant under rotation: $E_\mu(R\cdot\rho)=E_\mu(\rho)$ for $R \in \mathrm{SO}(3)$. 
Thus, equilibrium densities may form critical orbits, leading to rotation-induced degeneracy of the Hessian.

To compute the rotational directions, let
$e_1,e_2,e_3$ be the standard basis of
$\mathbb R^3$.
The corresponding orbit tangent vectors are
\begin{equation}
    z_\ell(\rho)(\mathbf p)
    =
    -\nabla_{\mathbb S^2}\rho(\mathbf p)
    \cdot(e_\ell\times\mathbf p),
    \qquad \ell=1,2,3,
\end{equation}
where $\nabla_{\mathbb S^2}$ denotes the surface gradient.
These vectors span the vertical space
$\mathcal V_\rho=T_\rho(\mathrm{SO}(3)\cdot\rho)$, 
which can be used to compute the horizontal projector for SR-HiSD.

\begin{figure}[H]
    \centering
    \includegraphics[width=\textwidth]{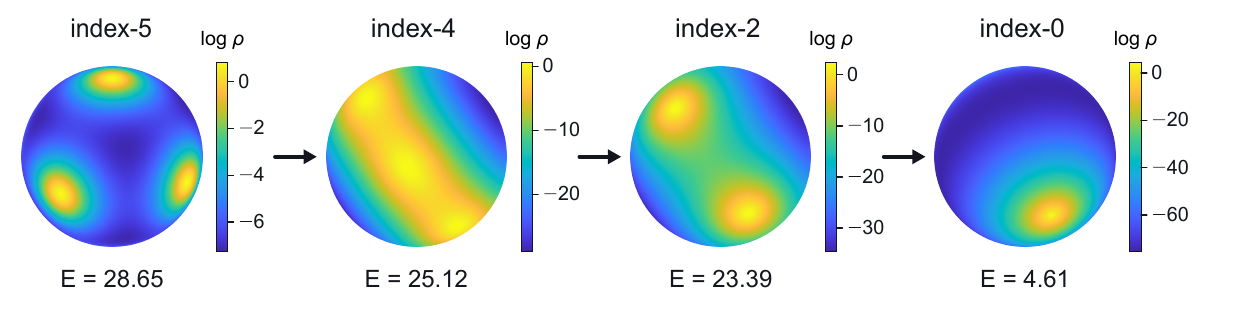}
    \caption{Downward search of the Onsager model at $\mu=80$ using SR-HiSD.}
    \label{fig:onsager_p80_downward_path}
\end{figure}

With this horizontal projector, we apply SR-HiSD to reproduce the downward search reported in \cite{yin2022solution}.
For $\mu=80$, we start from an index-5 saddle and obtain a sequence of critical points with indices $5\to4\to2\to0$, as shown in Figure~\ref{fig:onsager_p80_downward_path}.
The resulting pathway agrees with that reported in \cite{yin2022solution}, demonstrating that SR-HiSD recovers this downward search while excluding rotational directions from the saddle search.

\subsection{Landau--Brazovskii model with translational symmetry}
We finally consider the Landau--Brazovskii (LB) model on $\Omega=[0,L]^3$ with periodic boundary conditions \cite{zhang2026exploring}.
For model parameters $\tau$ and $\gamma$, the free energy is
\begin{equation}
    E(\phi)
    =
    \int_\Omega
    \left[
        \frac12|(\Delta+1)\phi|^2
        +\frac{\tau}{2}\phi^2
        -\frac{\gamma}{3!}\phi^3
        +\frac{1}{4!}\phi^4
    \right]\,\mathrm d\mathbf r.
\end{equation}
The order parameter $\phi$ satisfies the zero-mean constraint
$\int_\Omega\phi(\mathbf r)\,\mathrm d\mathbf r=0$, and belongs to the zero-mean periodic Sobolev space
$H=\left\{
        \phi\in H^2_{\rm per}(\Omega):
        \int_\Omega\phi(\mathbf r)\,\mathrm d\mathbf r=0
    \right\}$. 
Following \cite{zhang2026exploring}, we denote the HiSD method for this constrained LB model by LBSD.
We discretize the model using a Fourier pseudospectral method on a uniform $N^3$ grid with $N=128$.

The energy is invariant under periodic translations described by the group
$G=\mathbb T_L^3=\mathbb R^3/(L\mathbb Z)^3$.
Specifically, for the action $(a\cdot\phi)(\mathbf r)=\phi(\mathbf r-a)$,
we have $E(a\cdot\phi)=E(\phi)$.
At stationary states, the translational directions lie in the kernel of the Hessian and contribute to its degeneracy.
Away from stationary states, however, these directions need not lie in the Hessian kernel, unlike the rigid-translational directions \eqref{eqn:translation} in particle clusters.
To compute the horizontal projector, the orbit tangent
vectors are given by
$z_p(\phi)=\partial_{p}\phi$, $p=1,2,3$.
These vectors can be directly computed using Fourier methods and are already available from gradient calculations.

\begin{figure}[htbp]
    \centering
        \includegraphics[width=\linewidth]
        {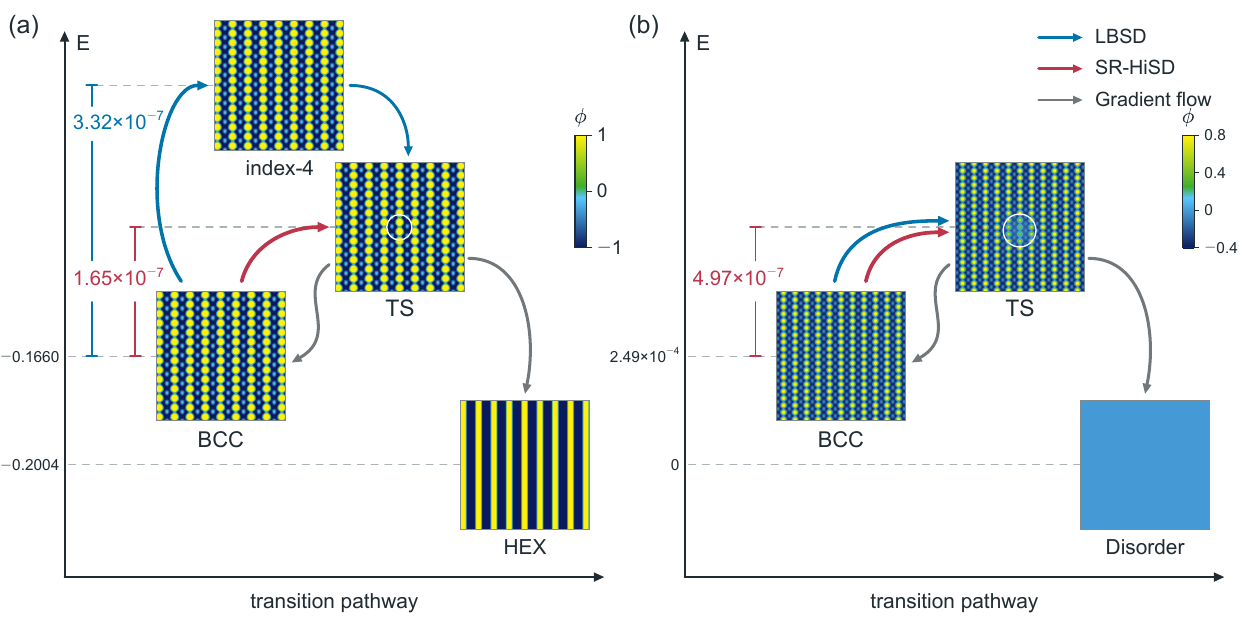}
    \caption{Comparison of LBSD and SR-HiSD for transitions (a) from BCC to HEX and (b) from BCC to disorder. Three-dimensional configurations are visualized through $L\times L$ square slices. TS denotes the transition state, and white outlines mark the critical nucleus. }
    \label{fig:LB_results}
\end{figure} 

First, we compute the nucleation pathway from the body-centered cubic state (BCC) to hexagonally-packed cylinder state (HEX) at $\tau=-0.3$, $\gamma=0.8$ and $L = 24\sqrt{2}\pi$. The Hessian at BCC has three zero eigenvalues associated with translations.
As shown in Figure~\ref{fig:LB_results}(a), 4-LBSD locates an index-4 saddle.
Subsequent downward searches using 1-LBSD then locate an index-1 saddle, which serves as the transition state (TS) with the critical nucleus. 
In contrast, 1-SR-HiSD directly locates the TS by restricting its ascent direction to the horizontal space. 
Gradient-flow trajectories initiated by perturbing the TS on both sides along its unstable eigenvector lead to BCC and HEX, thereby constructing the minimum energy pathway.

\begin{table}[htbp]
    \centering
    \caption{Computational cost comparison for the
    transitions in the LB model.}
    \label{tab:lb_cost_comparison}
    \begin{tabular}{lcccc}
        \hline
        Method & BCC to HEX & BCC to Disorder \\
        \hline
        4-LBSD    & 26,987 (15.68 h) & 35,696 (26.13 h) \\
        1-SR-HiSD & 15,372 (3.25 h) & 32,410 (6.76 h) \\
        \hline
    \end{tabular}
\end{table}

Next, we consider the transition from BCC to the disordered state $(\phi=0)$ at $\tau=0.05$, $\gamma=0.6$ and $L = 32\sqrt{2}\pi$.
Both 4-LBSD and 1-SR-HiSD can locate the same index-1 saddle connecting BCC and disorder, as shown in Figure~\ref{fig:LB_results}(b). 

Table~\ref{tab:lb_cost_comparison} compares the iteration counts and computational times of 4-LBSD and 1-SR-HiSD. 
Both methods employ the same BB step-size strategy and LOBPCG-based direction updates to ensure a fair comparison. 
Note that 4-LBSD converges to an index-4 saddle in the BCC-to-HEX transition. 
These results demonstrate the efficiency of SR-HiSD in locating saddle points of the prescribed Morse index.

\section{Conclusions}
\label{sec:5}
In this paper, we have proposed the SR-HiSD method for locating index-$k$ saddle points of $G$-invariant energy. 
Under suitable assumptions on the group action and Morse--Bott conditions, we formulate quotient HiSD on the quotient manifold.
For the quotient HiSD, we established local exponential stability of the target equilibrium states and instability of equilibrium states associated with critical points without a prescribed index. 
We then derive SR-HiSD in the ambient Hilbert space by taking the horizontal lift of the quotient HiSD, and prove that SR-HiSD locally converges to a saddle point with the prescribed index. 
We develop numerical algorithms for SR-HiSD by combining position updates with direction updates restricted to the horizontal space.
Numerical experiments on Morse clusters, the Onsager model, and the LB model demonstrated that SR-HiSD can locate saddle points of the prescribed Morse index and construct solution landscapes.
By excluding symmetry directions from saddle searches, SR-HiSD addresses the difficulties caused by symmetry-induced degeneracy and directly targets saddles of the prescribed index without introducing additional ascent directions.
The experiments also showed improved success rates in reaching the prescribed index and a substantial reduction in computational time in the large-scale LB example.
These results demonstrate the advantages of SR-HiSD in both saddle-point search and computational efficiency. 

The present analysis assumes a smooth quotient manifold under a free and proper group action. However, the isotropy group of a critical point may be nontrivial, and the quotient space may then have a stratified structure. Future work may extend SR-HiSD to this setting and study its local stability near such critical points. Discrete symmetries also deserve further study, particularly in chiral molecular systems, where symmetry-related states may be connected by physically meaningful transition paths. Extending SR-HiSD to account for these symmetries and identify the associated transition states is another possible direction for future research.

\bibliographystyle{siamplain}
\bibliography{references}
\end{document}